\documentclass[11pt,reqno]{amsart}
\usepackage{amsmath,amssymb,graphicx,amscd,amsfonts,psfrag,fancybox,qsymbols,indentfirst}
\usepackage{calligra}
\usepackage{mathrsfs}
\usepackage{dsfont}
\usepackage[latin1]{inputenc}

\DeclareMathAlphabet{\mathcalligra}{T1}{calligra}{m}{n}
 \DeclareMathOperator{\E}{\mathbb{E}}      
 \DeclareMathOperator{\Var}{\mathrm{Var}}      

 \newcommand{\ii}{{\mathrm{i}}}
  
\newcommand{\law}{\overset{\mbox{\rm \scriptsize law}}{=}}

\newcommand{\convlaw}{\overset{\mbox{\rm \scriptsize law}}{\longrightarrow}}

\newtheorem{thm}{Theorem}[section]

\newtheorem{lem}[thm]{Lemma}
\newtheorem{prop}[thm]{Proposition}
\newtheorem{rem}[thm]{Remark}
\theoremstyle{definition}

\theoremstyle{remark}

\def \be{\begin{eqnarray*}}
\def \ee{\end{eqnarray*}}
\def \ben{\begin{eqnarray}}
\def \een{\end{eqnarray}}

\def\sn{^{(n)}}

\numberwithin{equation}{section}

\begin{document}

\title[Bridges and random truncations of random matrices]{Bridges and random truncations of random matrices}
\date{\today}

\author{V. Beffara}
\address{UMPA -- ENS Lyon, UMR 5669, 
  46 allée d'Italie, 69364 Lyon cedex 07, France}
\email{vbeffara@ens-lyon.fr}

\author{C. Donati-Martin}
\address{Universit\'e  Versailles-Saint Quentin, LMV UMR 8100,
 B\^atiment Fermat, 45 avenue des Etats-Unis,
F-78035 Versailles Cedex}
\email{catherine.donati-martin@uvsq.fr}

\author{A. Rouault}
 \address{Universit\'e  Versailles-Saint Quentin, LMV UMR 8100,
 B\^atiment Fermat, 45 avenue des Etats-Unis,
F-78035 Versailles Cedex}
 \email{alain.rouault@uvsq.fr}


\maketitle

\begin{abstract}
  Let $U$ be a Haar distributed matrix in $\mathbb U(n)$ or $\mathbb O
  (n)$.  In a  previous paper,  we  proved that  after centering,  the
  two-parameter process
  \[T^{(n)} (s,t)  = \sum_{i \leq  \lfloor ns \rfloor, j  \leq \lfloor
    nt\rfloor}  |U_{ij}|^2\   ,  \   s,t  \in  [0,1]\]   converges  in
  distribution  to the  bivariate  tied-down Brownian  bridge. In  the
  present paper, we  replace the deterministic truncation of  $U$ by a
  random  one,  in which  each  row  (resp.\  column) is  chosen  with
  probability  $s$  (resp.\  $t$)  independently. We  prove  that  the
  corresponding   two-parameter    process,   after    centering   and
  normalization by $n^{-1/2}$ converges to  a Gaussian process. On the
  way we meet other interesting convergences.
\end{abstract}

\section{Introduction}

Let us consider a  unitary matrix $U$ of size $n\times  n$ and 
 two fixed integers  $p< n$  and $q < n$. Let us  call  $U^{p,q}$ the (rectangular)
matrix obtained by deleting the last  $n- p$ rows
and $n- q$  columns from $U$. It is well known that if  $U$ is Haar distributed
in $\mathbb  U(n)$, the random  matrix $U^{p,q}\left(U^{p,q}\right)^*$
has  a  Jacobi   matricial  distribution  and  that   if  $p,q$ and $n\rightarrow \infty$ 
 with $(p/n,  q/n)
\rightarrow  (s,t) \in(0,1)^2$,  its  empirical spectral  distribution
converges  to   a  limit   ${\mathcal  D}_{s,t}$  (see   for  instance
\cite{Collins}),   often    called   the    generalized   Kesten-McKay
distribution. 

In      \cite{CDMAR}      we       studied      the      trace      of
$U^{p,q}\left(U^{p,q}\right)^*$  which  is  also  the  square  of  the
Frobenius  (or  Euclidean)  norm  of $U^{p,q}$.  Actually  we  set  $p
=\lfloor ns\rfloor,  q= \lfloor nt\rfloor$ and  considered the process
indexed  by $s,t  \in [0,1]$.  We  proved that,  after centering,  but
without any  normalization, the process converges  in distribution, as
$n  \rightarrow \infty$,  to  a bivariate  tied-down Brownian  bridge.
Previously,  Chapuy   \cite{chapuy}  proved   a  similar   result  for
permutation matrices, with an $n^{-1/2}$ normalization.

Besides,  for purposes  of random  geometry analysis,  Farrell  has
proposed    another    model    in   \cite{Farfourier}    (see    also
\cite{Fararxiv}), deleting randomly and independently a proportion $1-
s$ of rows and a proportion $1-  t$ of columns from a Haar distributed
matrix in $\mathbb U (n)$. If $\mathcal U^{s,t}$ denotes the matrix so
obtained, 
he proved that (for fixed $s,t$) the empirical spectral
distribution of $\mathcal  U^{s,t}\left(\mathcal U^{s,t}\right)^*$ converges again, as $n \rightarrow \infty$, to  ${\mathcal D}_{s,t}$.

It is then  tempting to study the
trace of $\mathcal  U^{s,t}\left(\mathcal U^{s,t}\right)^*$ as a process, after having defined 
 a probability space where  all  random truncations live  simultaneously. 
 For that purpose we define a  double array of $2n$ auxiliary independent uniform variables $R_1, \dots, R_n, C_1, \dots,
 C_n$ and then,  for any choice of $(s,t)$, obtain the matrix $\mathcal U^{s,t}$ by removing from $U$
 rows with indices not in $\{i: R_i \leq s\}$ and
 columns with indices not in $\{j : R_j \leq t\}$. This gives us a coupled realization of the $\mathcal U^{s,t}$, reminiscent of the ``standard coupling'' for percolation models. Then, we notice  that in
 the first model, the invariance  of the Haar distribution on $\mathbb
 U (n)$  implies that we  could have deleted any  fixed set of  $n- p$
 rows and  $n-q$ columns.  So, we can  consider the  random truncation
 model  as  the  result  of the  subordination  of  the  deterministic
 truncation model by  a couple of binomial processes.  In other words,
 we treat the latter uniform variables as an environment, and
 state 
 quenched and annealed convergences.

For instance, we will prove that after convenient centering and without normalization, the above process converges (quenched) to a bivariate Brownian bridge, but that after another centering and with normalization by $n ^{-1/2}$ it converges (annealed) to 
 a Gaussian 
 process which is no more a bivariate Brownian bridge.

 We use the 
 space
$D([0,1]^2)$  endowed with the
topology  of  Skorokhod   (see  \cite{bickel1971convergence}). It consists of functions from
$[0, 1]^2$  to $\mathbb R$  which are  at each point  right continuous
(with respect  to the natural partial  order of $[0, 1]^2$)  and admit
limits in all ``orthants''. 
 For the sake of completeness, we treat also the one-parameter process,
i.e.\ truncation  of the first column  of the unitary matrix,  and the
case of permutation matrices.

Actually,  Farrell  considered first  the
(deterministic) discrete Fourier transform (DFT) matrix
\begin{equation}\label{DFTdef}F_{jk} = \frac{1}{\sqrt n} e^{-2\ii  \pi (j-1)(k-1)/n} \ , \quad j,k
= 1, \dots, n\,,\end{equation} and proved that  after random truncation, a Haar unitary matrix has the same limiting singular value distribution. In  a still more  recent paper (\cite{AndFar}),  Anderson and Farrell
explain the connection with liberating sequences. In some sense, the randomness coming from the truncation is stronger than the randomness of the initial matrix. Here, we have considered also
 the DFT  matrix, but we can
as well consider any (random or  not random) matrix whose elements are
all  of  modulus  $n^{-1/2}$,  for  instance  a  (normalized)  complex
Hadamard matrix.

The paper is organized as follows. In Sec. 2 we provide some definitions. Section 3 is devoted to our main results, the convergence of one-parameter (Theorem~\ref{onedim}) and two-parameter processes (Theorems~\ref{q-a} and \ref{twodimpermut}). In Sec. 4, we introduce the subordination method, which allows to give  the proofs of the latter theorems as examples of application.  
  In  Section   5,  we  go  back   to  the  direct  method,   used  in
  (\cite{arXiv}) which does not assume that the result of deterministic truncation is known. This point of view leads to conjectures.

\section{Notation}

We introduce the random processes that  we will consider in this paper
and the various limiting processes involved.

\subsection{The main statistics}
Let $\mathbb  U(n)$ (resp. $\mathbb O(n)$) be the  group of unitary (resp. orthogonal)   $n\times n$ matrices and  $U = (U_{ij})$ its generic element.
 We equip $\mathbb U(n)$  (resp. $\mathbb O(n)$) with the
Haar probability measure $\pi\sn(dU)$.

\noindent To define
two systems  of projective Bernoulli  choices of rows and  columns, we
will  need two  independent families  of independent  random variables
uniformly distributed on $[0,1]$ so  that we can treat the randomness
coming from the  truncation as an environment.  More specifically, the
space of environments is $\Omega= [0,1]^\mathbb N \times [0,1]^\mathbb
N$, whose generic element  is denoted by $\omega = (R_i, i  \geq 1, C_j, j
\geq 1)$. We equip $\Omega$ with the probability  measure $d\omega$ which is  the infinite product
of copies of the uniform distribution. In the sequel, ``for almost every $\omega$'' will mean
``for $d\omega$ - almost every $\omega$.''

For the one-parameter model, we introduce two processes with values in
$D([0,1])$:
\begin{eqnarray}
  B\sn_s (U)&=& \sum_1^{\lfloor ns\rfloor} |U_{i1}|^2  \,,\\
  \label{repres1}
  {\mathcal B}_s\sn(\omega, U)  &=& \sum_1^n |U_{i1}|^2 1_{R_i \leq s}\,.
\end{eqnarray}

For the  two-parameter model,  we introduce  processes with  values in
$D([0,1]^2)$:
\begin{enumerate}
\item $T\sn(U)$
  defined by
  \[T\sn_{s,t}(U) =  \sum_{i=1}^{\lfloor ns\rfloor}\sum_{j=1}^{\lfloor
    nt\rfloor} |U_{ij}|^2\,,\]
\item $\mathcal T\sn(\omega, U)$ defined by \ben
  \label{defcalT}
  \mathcal    T\sn_{s,t}(\omega,    U)   =    \sum_{i=1}^n\sum_{j=1}^n
  |U_{ij}|^2\mathbf 1_{R_i \leq s} \mathbf 1_{C_j \leq t}\,.\een
\end{enumerate}
 The counting  processes  $S\sn$ and  ${S'}\sn$  are defined  by
  \ben\label{Sandco} S\sn_s(\omega) = \sum_{i=1}^n \mathbf 1_{R_i \leq
    s} &,& {S'_t}\sn(\omega) = \sum_{j=1}^n \mathbf 1_{C_j \leq t}\,,\een
  and    their    normalized    version    $\widetilde    S\sn$    and
  $\widetilde{S'}\sn$ by \ben
 \label{Sandcon}
 \widetilde   S\sn_s   =    n^{-1/2}\left(S\sn_s   -   ns\right)   &,&
 \widetilde{S'}_t\sn= n^{-1/2}\left({S'}_t\sn - nt\right)\,. \een

\subsection{Gaussian processes and bridges}

The classical Brownian bridge denoted  by $B_0$ is a centered Gaussian
process with continuous paths defined on $[0,1]$, of covariance
\[\mathbb E \left(B_0(s)B_0(s')\right) =  s\wedge s' - ss'\,.\]
The  bivariate Brownian  bridge  denoted by  $B_{0,0}$  is a  centered
Gaussian  process  with  continuous  paths  defined  on  $[0,1]^2$  of
covariance
\[\mathbb   E   \left(B_{0,0}(s,t)B_{0,0}(s',t')\right)  =   (s\wedge
s')(t\wedge t')  - ss'tt'.\]  The tied-down bivariate  Brownian bridge
denoted  by  $W^{(\infty)}$  is   a  centered  Gaussian  process  with
continuous paths defined on $[0,1]^2$ of covariance
\[\mathbb  E [W^{(\infty)}(s,t)  W^{(\infty)}(s',t')]  = (s\wedge  s'-
ss')(t\wedge t'  - tt').\] Let  also ${\mathcal W}^{(\infty)}$  be the
centered Gaussian  process with continuous paths  defined on $[0,1]^2$
of covariance
\[\mathbb      E      [{\mathcal     W}^{(\infty)}(s,t)      {\mathcal
  W}^{(\infty)}(s',t')] =ss'(t\wedge t') + (s\wedge s')tt' -2ss'tt'.\]
It can be defined also as \ben
\label{defwinfini}{\mathcal   W}^{(\infty)}(s,t)   =   sB_0(t)   +   t
B'_0(s)\een  where  $B_0$ and  $B'_0$  two  independent one-parameter
Brownian bridges.

At last we will meet the  process denoted by $B_0\otimes B_0$ which is
a centered process with continuous paths defined on $[0,1]^2$ by
\[B_0 \otimes B_0 (s,t) = B_0(s)  B'_0(t)\] where $B_0$ and $B'_0$ are
two independent Brownian bridges. This process is not Gaussian, but it
has the same covariance as $W^{(\infty)}$.

Similarly, if $F$ and $G$ are two processes with values in $D([0,1])$, defined on the same probability space, we denote by
$F\otimes G$ the process with values in $D([0,1]^2)$ defined by
\[F\otimes G(s,t) = F(s) G(t).\]
For simplicity we denote by
$I$ the deterministic trivial process $I_s=s$.

\section{Convergence in distribution}

We present unified results in the  cases of the unitary and orthogonal
groups. For this purpose we use the classical notation
\[\beta' = \frac{\beta}{2}= \begin{cases}1/2 \;\;\;\;&\mbox{in the orthogonal case} ,\\
  1\;\;&\mbox{in the unitary case.}
\end{cases}\] 

\subsection{One-parameter processes}

Let us begin with the  one-parameter processes, where $\convlaw$ means
convergence in distribution in $D([0,1))$. We present successively the results for the deterministic and random truncations.
\begin{lem}
  \label{silver}
  Under $\pi\sn$,
  \begin{eqnarray}
    \label{pont1a}
    n^{1/2}\left(B\sn - I\right)\convlaw \sqrt{\beta'^{-1}}B_0\,.
  \end{eqnarray}
\end{lem}

This convergence \eqref{pont1a} is well known since at least Silverstein \cite{Silver1}
(in the case $\beta'=1$). It can  be viewed as a direct consequence of
the fact  that the  vector $(|U_{i1}|^2, i=1,  \dots, n)$  follows the
Dirichlet $(\beta', \dots, \beta')$ distribution on the simplex.

\begin{thm}
  \label{onedim}
  \begin{enumerate}
  \item  (Quenched)  For almost  every  $\omega$,  the push-forward  of
    $\pi\sn(dU)$ by the map \ben
    \label{qone}
    U  \mapsto  n^{1/2}   \left(\mathcal  B\sn(\omega,  U)  -  n^{-1}
      S\sn(\omega)\right) \een converges weakly to the distribution of
    $\sqrt{\beta'^{-1}}B_0$.
  \item   (Annealed)   Under   the  joint   probability   measure $d\omega\otimes
    \pi\sn(dU)$
    \begin{eqnarray}
      \label{pont1b}
      n^{1/2} \left({\mathcal B}\sn - I\right)\convlaw \sqrt{1+ \beta'^{-1}} B_0.
    \end{eqnarray}
  \end{enumerate}
\end{thm}

\subsection{Two-parameter processes}

Let us continue with the two-parameter processes, where now $\convlaw$
means convergence in distribution in $D([0,1]^2)$. We
study  three models. In the first one, $U$ is   the DFT matrix defined in (\ref{DFTdef}). In the second one, $U$ is sampled from the Haar measure on $\mathbb U(n)$ or $\mathbb O(n)$. Though the proof is much more involved than in the first model, the annealed convergence gives the same limit. At last, 
for the sake  of completeness, we state here a result when 
$U$ is chosen uniformly among $n \times n$ permutation matrices.

\subsubsection{DFT}

Here,  there is  no  randomness in  $U$, so  that  we have the decomposition:  
\ben\label{trivialmathcalT}
\mathcal T\sn=
\frac{S\sn\otimes{S'}\sn}{n}= n  I\otimes I + \widetilde  S\sn \otimes
\widetilde{S'}\sn + n^{1/2}\left(\widetilde S\sn  \otimes I + I\otimes
  \widetilde{S'}\sn\right)\!. \een

\begin{thm}
  \label{twodimDFT}
  If $U$  is the DFT  matrix (or  more generally  if $U$ is any matrix such that  $|U_{ij}|^2 =1/n$
  a.s.\ for every $i,j$), then under the probability measure  $d\omega$
\ben
  \label{convDFT}n^{-1/2}\left(\mathcal  T\sn  -  \mathbb  E  \mathcal
    T\sn\right) \convlaw {\mathcal W}^{(\infty)}\,.\een
\end{thm}

\proof It  is straightforward  since the processes  $S\sn$ and  ${S'}\sn$ are
independent. From (\ref{trivialmathcalT}) we have successively
\[\mathbb      E     \mathcal      T\sn      =     n^{-1}      \mathbb
E\left(S\sn\otimes{S'}\sn\right) = n I\otimes I,\]
\[n^{-1/2}\left(\mathcal  T\sn  -  \mathbb E  \mathcal  T\sn\right)  =
n^{-1/2}  \left(\widetilde S\sn  \otimes \widetilde{S'}\sn\right) +\widetilde
  S\sn  \otimes  I  + I\otimes  \widetilde{S'}\sn\,.\]
Applying Donsker's theorem, we get
\[( \widetilde S\sn  , \widetilde{S'}\sn)\convlaw (B_0, B'_0)\]
so that
\[n^{-1/2} \left(\widetilde S\sn  \otimes \widetilde{S'}\sn\right) \rightarrow 0\]
in probability and
\[ S\sn  \otimes  I  + I\otimes  \widetilde{S'}\sn \convlaw B_0\otimes I + I \otimes B'_0 = \mathcal
W^{(\infty)}\,.\]
 $\Box$

\subsubsection{Haar unitary or orthogonal matrices}

 The
case of deterministic truncation was treated in our previous paper and
recalled now.
\begin{thm}[\cite{CDMAR}]
  \label{cdmar}
  Under $\pi\sn$, \ben
  \label{rappelnous}
  W\sn   :=  T\sn   -\mathbb   E   T\sn  \convlaw   \sqrt{\beta'^{-1}}
  W^{(\infty)}\,.\een
\end{thm}

\noindent The case of random truncation is ruled by the following theorem, which is the main result of the present paper.
\begin{thm}
  \label{q-a}
  \begin{enumerate}
  \item  (Quenched)  For almost  every  $\omega$,  the push-forward  of
    $\pi\sn(dU)$ on $D([0,1])^2$ by the mapping \ben
    \label{defcalV}U \mapsto \mathcal V\sn := \mathcal T\sn(\omega, U)
    -  \frac{S\sn(\omega)  \otimes  {S'}\sn(\omega)}{n}\een  converges
    weakly to the distribution of $\sqrt{\beta'^{-1}} W^{(\infty)}$.
  \item (Annealed)  Under the joint probability measure 
 $d\omega \otimes\pi\sn(dU)$, 
  \[ n^{-1/2} \left(\mathcal T\sn  - \mathbb E \mathcal T\sn
\right)
 \convlaw \mathcal W^{(\infty)}\,.\]
  \end{enumerate}
\end{thm}
\begin{rem}
  \label{remn0}
  Let $M^{p,q} =  U^{p,q}\left(U^{p,q}\right)^*$ and $\mathcal M^{s,t}
  = \mathcal U^{s,t}\left(\mathcal U^{s,t}\right)^*$. For $s,t$ fixed,
  the  random variables  $T\sn_{s,t}$  and  $\mathcal T\sn_{s,t}$  are
  linear  functionals  of  the   empirical  spectral  distribution  of
  $M^{\lfloor  ns\rfloor, \lfloor  nt\rfloor}$ and  $\mathcal M^{s,t}$
  respectively.  For classical  models  in Random  Matrix Theory,  the
  convergence of fluctuations of such linear functionals do not need a
  normalizing factor,  since the variance is  bounded (the eigenvalues
  are  repelling  each other).  Here,  this  is  indeed the  case  for
  $T\sn_{s,t}$ (see \cite{DuPa} for  the complete behavior for general
  tests functions). But, in the case of $\mathcal T\sn_{s,t}$, we have
  $  \Var \left(\mathbb  E [\mathcal  T\sn_{s,t} |  \omega ]\right)  =
  O(n)$, which demands a  normalization. Notice however that the
    main  source of  this variance  lies  in the  fluctuations of  the
    number  of columns  and  lines removed  from  the initial  matrix,
    rather than in the matrix itself.
\end{rem}
\subsubsection{Permutation matrices}

 Let us call $p\sn$
the  uniform  measure  on  the group  $\mathcal  S_n$  of  permutation
matrices of $\{1, \dots,n\}$. The deterministic truncation was treated
by Chapuy.
\begin{thm}[\cite{chapuy}]
  \label{permut}
  Under $p\sn$ we have
  \begin{eqnarray}
    \label{chap}n^{-1/2}\left(T\sn  -\mathbb  E  T\sn\right)&\convlaw&
    W^{(\infty)}.
  \end{eqnarray}
\end{thm}

Here  is  the  result  for  the  statistics  obtained  by  the  random
truncation.
\begin{thm}
  \label{twodimpermut}
  \begin{enumerate}
  \item  (Quenched)  For almost  every  $\omega$,  the push-forward  of
    $p\sn$ 
    by the mapping
    \[U    \mapsto    n^{-1/2}\left(\mathcal   T\sn(\omega,    U)    -
      \frac{S\sn(\omega)      \otimes     {S'}\sn(\omega)}{n}\right)\]
    converges weakly to $W^{(\infty)}$.
  \item (Annealed) Under the  joint probability measure $d\omega\otimes p\sn(dU)$
    \ben\label{nchap1}n^{-1/2}\left(\mathcal T\sn - \mathbb E \mathcal
      T\sn\right) \convlaw B_{00}\,.\een
  \end{enumerate}
\end{thm}

\section{Proofs by subordination}

We   present    here   proofs    of   Theorems~\ref{onedim},~\ref{q-a}
and~\ref{twodimpermut} whose 
key point is 
a representation by subordination.

\subsection{Preliminaries}

\begin{prop}
  \label{bef1}
  Assume that  $U$ is  a random  unitary matrix  such that  the matrix
  whose generic entry is $|U_{ij}|^2$  has a distribution invariant by
  multiplication (right or left)  by permutation matrix. Let $\widehat
  T\sn$ be defined by
\[\widehat  T\sn_{s,t}(\omega,  U)  = T\sn_{n^{-1}S\sn_s  (\omega),
    n^{-1}{S'_l}\sn  (\omega)}(U)\,.\]
  Then  for  every  $\omega$  the
  push-forward of $\pi\sn(dU)$ by the  mapping $U \mapsto \mathcal T\sn
  (\omega, U)$ is  the same as the push-forward of  $\pi\sn(dU)$ by the
  mapping $U \mapsto \widehat T\sn(\omega, U)$. As a result the law of
  $\mathcal T\sn$ and $\widehat T\sn$ have the same distribution under
  $ d\omega\otimes\pi\sn(dU)$.
\end{prop}

\proof

Let  $R=  (R_1, \dots,  R_n)$  and  $C =  (C_1,  \dots,  C_n)$ be  two
independent samples of uniform variables on $[0,1]$. The corresponding
reordered samples are $\widetilde R  = (R_{(1)} , \dots, R_{(n)})$ and
$\widetilde C = (C_{(1)}, \dots,  C_{(n)})$, and the associated random
permutations are $\sigma$ and $\tau$, are defined by
\[R_{(i)} = R_{\sigma^{-1}(i)} \ , \ C_{(j)} = C_{\tau^{-1}(j)}, \ i,j
= 1, \dots, n\,.\] Moreover  $\sigma$ and $\widetilde R$ (resp. $\tau$
and $\widetilde C$) are independent. With these notations, we have
\begin{eqnarray*}{\mathcal T}\sn_{s,t} &=& \sum_{i,j=1}^n |U_{ij}|^2 1_{R_i \leq s} 1_{C_j \leq t}=  \sum_{i,j=1}^n
  |U_{\sigma^{-1}(i) \tau^{-1}(j)}|^2 1_{R_{(i)}\leq s} 1_{C_{(j)}\leq t}\\
  &=& \sum_{i \leq S\sn_s , j \leq {S'}\sn_t}  |U_{\sigma^{-1}(i) \tau^{-1}(j)}|^2 = T\sn_{n^{-1}S_s\sn, n^{-1}{S'}_t\sn}(\sigma U\tau^{-1})\,,
\end{eqnarray*}
where we  have identified  the permutations  $\sigma$ and  $\tau$ with
their matrices. Let $F$  be some test  function from  $D([0,1]^2)$ to $\mathbb  R$. We
have
\[\mathbb E [F(\mathcal T\sn(\omega, U))|\omega] = \mathbb E[F(T\sn_{n^{-1}S\sn, n^{-1}{S'}\sn}(\sigma U\tau^{-1})) | \omega]\,.\]
Since
the distribution of $(|U_{ij}|^2)_{i,j= 1}^n$ is invariant by permutation we get
\[\mathbb E [F(\mathcal T\sn(\omega, U)|\omega] =  \mathbb E[F(T\sn_{n^{-1}S\sn, n^{-1}{S'}\sn}(U)) | \omega]\]
or, in other words
\[\mathbb E [F(\mathcal T\sn(\omega, U)|\omega] = \mathbb E [F(\widehat T\sn(\omega, U)|\omega]\,,\]
which ends the proof. $\Box$

Now, the  key point to  manage the  subordination of processes  is the
following proposition.
\begin{prop}
  \label{subordination}
  Let $d$ be $1$ or $2$ and let $A\sn$ be a sequence of processes with
  values in $D([0,1]^d)$  such that $A\sn \convlaw A$.  Let $S\sn$ and
  ${S'}\sn$ be two independent  processes defined as in (\ref{Sandco})
  and independent of $A\sn$.
\begin{itemize}
\item If  $d=1$, set  $\mathcal A\sn  = \left(A\sn\left(n^{-1}S\sn_s)
      \right)\ ,\ s \in [0,1]\right)$. Then
    \[\left(\mathcal  A\sn  ,  \widetilde  S\sn\right)  \convlaw  (A  ,
    B_0);\]
 \item If  $d=2$, set $\mathcal A\sn  = \left(A\sn\left(n^{-1}S\sn_s,
        n^{-1} {S'_t}\sn \right)\ ,\ s, t \in [0,1]\right)$. Then
    \[\left(\mathcal A\sn  , \widetilde  S\sn ,  \widetilde {S'}\sn\right)
    \convlaw (A , B_0, B'_0),\]
  \end{itemize}
  where $A,  B_0, B'_0$ are independent  and $B_0$ and $B'_0$  are two
  independent (one-parameter) Brownian bridges.
\end{prop}
Notice  that  the marginal  convergence  of  $\widetilde S\sn$  (or  of
$(\widetilde S\sn, \widetilde {S'}\sn)$) is nothing but Donsker's theorem.

\proof

Let us  restrict us to  the case $d=1$  for simplicity. We  follow the
lines of proof of Theorem 1.6 of Wu \cite{Biao}. According to the
Skorokhod representation theorem, we can build a probability space and
stochastic processes $\mathbf  A\sn$, $\mathbf A$, 
$\widetilde{\mathbf S}\sn$, $\mathbf B_0$ on it such that
\begin{itemize}
\item all processes are $D([0,1])$ valued
\item $\mathbf A\sn$ and 
$\widetilde{\mathbf S}\sn$ are independent and
  $\mathbf   A\sn  \law   A\sn$   ,  
$\widetilde{\mathbf S}\sn \law \widetilde S\sn$
\item $\mathbf  A$ and  $\mathbf B_0$ are  independent and  $\mathbf A
  \law A$ , $\mathbf B_0 \law B_0$
\item $\mathbf  A\sn$ and 
$\widetilde{\mathbf S}\sn$  converge a.s.\ to
  $\mathbf A$ and $\mathbf B_0$, respectively.
\end{itemize}

Set
${\mathbf S}\sn = n^{1/2}\widetilde{\mathbf S}\sn + n I$. The  convergence  a.s.\  of
  $\widetilde{\mathbf S}\sn$ entails  the
convergence   a.s.\   of   $n^{-1}\mathbf   S\sn$  to   $I$   as   a
$D([0,1])$-valued   non-decreasing   process.    Here   the   limiting
subordinator  is  continuous. We  are  exactly  in the  conditions  of
Theorem 1.2 of \cite{Biao}, which allows to say
\[a.s.-  \lim_n \mathbf A_n \circ \mathbf S\sn = \mathbf A \]
and of course
\[a.s. -  \lim_n\ \! (\mathbf A_n  \circ \mathbf S\sn  , \widetilde{\mathbf
  S}\sn) = (\mathbf  A, \mathbf B_0)\,.\] Now, we  conclude, going down
to the convergence in distribution,
\[(\mathcal A_n  , \widetilde{S}\sn)  \convlaw (A, B_0)\,,\]  where $A$
and $B_0$ are independent. $\Box$
\subsection{Proof of Theorem~\ref{onedim}}

From Proposition~\ref{bef1} (stated for one-parameter  processes), we
have the equality  in law (conditionally on  $\omega$) \ben \{\mathcal
B_s\sn(\omega,.)\,, s \in [0,1]\}\law \{B\sn_{n^{-1}S_s\sn(\omega)}(.)
, s\in [0,1]\} \een
and                 then                  we                 decompose
\ben\label{decomp1n}n^{1/2}\left(B\sn_{n^{-1}S_s\sn(\omega)} -s\right)
=   n^{1/2}\left(B\sn_{n^{-1}S_s\sn(\omega)}  -n^{-1}S\sn_s\right)   + \widetilde S\sn_s\,.
 \een
If we set $A\sn(s) = n^{1/2} \left(B\sn_s
  -n^{-1}\lfloor ns \rfloor\right)$,
Lemma~\ref{silver} above  says that we are exactly  in the assumptions
of   Proposition~\ref{subordination}.  Both   processes  of   the  RHS
of~(\ref{decomp1n}) converge  in distribution towards  two independent
processes,   distributed   as    $\sqrt{\beta'^{-1}}B_0$   and   $B_0$
respectively,   hence   the   sum   converges   in   distribution   to
$\sqrt{1+\beta'^{-1}}B_0$. $\Box$

\subsection{Proofs of Theorems~\ref{q-a} and~\ref{twodimpermut}}

From Proposition~\ref{bef1},  we reduce the  problems to the  study of
$\widehat T\sn$.
Let us first remark that \ben
\label{econd}
\mathbb  E [\mathcal  T\sn  | \omega]  = \mathbb  E  [\widehat T\sn  |
\omega] =  n^{-1} S\sn(\omega)  \otimes {S'}\sn(\omega)\,. \een  If we
set \ben\label{defY} \widehat W\sn (\omega, U) = \widehat T\sn(\omega,
U)  -   \mathbb  E  [\widehat  T\sn   |  \omega]  \een  we   have  the
decomposition: \ben
\label{decnew}
\widehat T\sn -  \mathbb E \widehat T\sn = \widehat  W\sn + \widetilde
S\sn     \otimes     \widetilde{S'}\sn     +     n^{1/2}\left(I\otimes
  \widetilde{S'}\sn + \widetilde S\sn \otimes I\right)\,. \een

\subsubsection{Proof of Theorem~\ref{q-a}}
For the  quenched fluctuations of  $\widehat W\sn$, we are  exactly in
the  assumptions  of  Proposition~\ref{subordination},  with  $A\sn  =
\widehat  W\sn$  and  $A= \sqrt{\beta'^{-1}}W^{(\infty)}$,  thanks  to
Theorem~\ref{cdmar}.
This implies in particular that {\it (1)\/} holds.

For {\it (2)}, from  Proposition~\ref{subordination}, we see also that
\ben
\label{vb3}
(\widehat   W\sn,   \widetilde   S\sn,   \widetilde{S'}\sn)   \convlaw
(\sqrt{\beta'^{-1}}W^{(\infty)},  B_0, B'_0)\,  \een  where the  three
processes are independent. This implies
\[ \widehat W\sn +  \widetilde S\sn \otimes \widetilde{S'}\sn \convlaw
\sqrt{\beta'^{-1}}W^{(\infty)} + B_0\otimes B'_0\,\] and consequently
\[n^{-1/2}\left(   \widehat    W\sn   +   \widetilde    S\sn   \otimes
  \widetilde{S'}\sn\right) \rightarrow 0\]  in probability. Looking at
the decomposition (\ref{decnew}) and using again the convergence
\[I\otimes  \widetilde{S'}\sn +  \widetilde  S\sn  \otimes I  \convlaw
\mathcal W^{(\infty)}\] we conclude that
\[n^{-1/2}\left(\widehat  T\sn   -  \mathbb  E   \widehat  T\sn\right)
\convlaw \mathcal W^{(\infty)}\] which  is equivalent to the statement
of {\it (2)}. $\Box$

\subsubsection{Proof of Theorem~\ref{twodimpermut}}

We are now in the assumptions of Proposition~\ref{subordination}, with
$A\sn =  n^{-1/2}(T\sn -  \E(T\sn))$ and  $A= W^{(\infty)}$  thanks to
Theorem~\ref{permut}.  This implies  in  particular  that {\it  (1)\/}
holds.

For {\it (2)}, from  Proposition~\ref{subordination}, we see also that
\ben
\label{vb3p}
\left(n^{-1/2}\widehat    W\sn,   \widetilde    S\sn,   \widetilde{S'}
  \sn\right) \convlaw (W^{(\infty)}, B_0, B'_0)\, \een where the three
processes are independent. This implies
\[n^{-1/2}\left(\widehat  T\sn   -  \mathbb  E   \widehat  T\sn\right)
\convlaw  W^{(\infty)}  +  \mathcal W^{(\infty)}\,,\]  where  the  two
processes in the RHS are independent. The equality in law
\[B_{0,0}  \law  W^{\infty}  +  \mathcal W^{\infty}\]  was  quoted  in
\cite{DPY} section 2. $\Box$

\begin{rem}
In the way leading from Theorem~\ref{cdmar} to Theorem~\ref{q-a},
we can  see that the Haar  distribution of random unitary  matrices is
not  involved in  the  proofs, except  via the  invariance  in law  by
permutation   of   rows   and   columns  (which   is   the   core   of
Proposition~\ref{bef1}). In a  recent work, Benaych-Georges \cite{FBG}
proved that, under some conditions, the unitary matrix of eigenvectors
of a  Wigner matrix induces the  same behavior for the  asymptotics of
$T\sn$ as  under the Haar  distribution. In the same  vein, Bouferroum
\cite{AliB}  proved a  similar  statement for  the  unitary matrix  of
eigenvectors of  a sample  covariance matrix.  We could  ask if  it is
possible to take benefit of these  results to give an extension of our
theorem  to  more general  random  unitary  matrices. Actually,  these
authors proved that if the eigenvalues are ordered increasingly and if
we call  $U_<$ the  matrix of  corresponding eigenvectors,  under some
assumptions, the process $T\sn (U_<)-  \mathbb E T\sn (U_<)$ converges
in law to $\sqrt{\beta'^{-1}}W^{(\infty)}$ as in Theorem~\ref{cdmar}. To be able to apply
Proposition~\ref{bef1}, we would have to check  the invariance of the law of
the matrix $(|(U_<)_{ij}|^2)$ under  multiplication (left or right) by
a  permutation matrix.  A  short look  reveals 
 that  if  $\sigma$ is  a
permutation,  then in  both  models $\sigma  M\sn\sigma^*$ and  $M\sn$
share  the  same  eigenvalues,  and  $\sigma U_<$  is  the  matrix  of
eigenvalues of $\sigma  M\sn\sigma^*$. But they may not  have the same
distribution if they  are complex. Even if we restrict  to real Wigner
matrices, we have  indeed $\sigma U_< \law U_<$ but  the other type of
equality  in law  (right  permutation)  is in  general  not true.  For
further remarks on  the type of unitary matrices which  could give the
same convergence, see Section \ref{about}.
\end{rem}

\section{About direct proofs of the main results and two conjectures}
\label{about}
First, let us remark that in  the Haar and permutation models, the key
tool of the above approach  to random truncation was the subordination
machine.  It  assumes  that  we  know  the  previous  results  on  the
deterministic  truncation. Going  back  to the  proof  of this  latter
result  in the  Haar case \cite{CDMAR},  we see  that estimates  of moments  of all
degrees of monomials  in entries of the unitary matrix  are needed. We
can ask if a direct method  to tackle the random truncation demands so
high  moments  estimates.  This  fact, among  others,  legitimates  an
interest for direct proofs, starting from the
representation
(\ref{defcalT})  in  the Haar  and  permutation  models and  from  the
representation (\ref{repres1}) for the one-dimensional process.

A second  striking fact in the  study of the two-parameter  process is
that in 
Theorems~\ref{twodimDFT} and~\ref{q-a}(2), the  limiting processes are
the same. In other words, the behavior of the sequence of DFT matrices
is the same as  the mean behavior of sequence of  Haar matrices. If we
define $\mathbb U(\infty) := \times_{n=1}^\infty \mathbb U(n)$, we can
then ask how large is the set
\[\mathcal E := \{u =  (U\sn,  n \in \mathbb N) \in \mathbb U(\infty)
\  | \  n^{-1/2}(\mathcal T\sn  (.  , U\sn)  - n  I\otimes I)\convlaw
\mathcal W^{(\infty)}\}\,.\]
Actually, it is  equivalent to consider $\omega$ as  the random object
and the collection of $\mathbb U(n), n  \in \mathbb N$ as the space of
environments.
There  are  several  choices  to  equip  $\mathbb  U(\infty)$  with  a
probability measure  such that its  marginal on $\mathbb U(n)$  is the
Haar measure $\pi\sn$.
In \cite{B2N}, the authors introduced  the notion of virtual isometry.
They consider  a family of projections  ${\mathcalligra p}_{m,n}$ from
$\mathbb U\sn$ to  $\mathbb U^{(m)}$ for $m <n$ and  define the subset
$\mathcal  U(\infty)$   of  $\mathbb  U(\infty)$  of   $u$  such  that
$\mathcalligra p_{m,n} (U\sn) = U^{(m)}$ for every $m,n$ with $m < n$.
They conclude that there exists  a unique probability measure $\pi$ on
$\mathcal U(\infty)$ (equipped  with the cylindrical $\sigma$-algebra)
whose $n$-th marginal is $\pi\sn$ for every $n$. Their construction is
also  compatible with  the framework of permutations  
(replace  $\pi\sn$ by
$p\sn$ and  $\pi$ by $p$). It  could be noticed that  the DFT sequence
belongs to $\mathbb U(\infty) \setminus \mathcal U(\infty)$. Besides,  in  \cite{JiangMax},  Jiang    ``inspired  by  a  common
statistical procedure  for simulating  a sequence of  Haar distributed
matrices  in statistical  programs''   assume that  $(U\sn  , n  \in
\mathbb N)$ is an independent sequence. The same remarks hold for $\mathbb O(\infty):= \times_{n=1}^\infty \mathbb O(n)$. 

In  the  following  subsections,  we will  give  alternate  proofs  of
Theorems~\ref{onedim}   and~\ref{q-a}  (annealed). We  propose also  two
conjectures  about weak  convergences  
 conditionally  upon $u$, for $\pi$ - almost every $u  \in \mathbb U(\infty)$ (resp. $\mathbb O(\infty)$), where $\pi$ is   any probability measure  whose
marginals are Haar measures on  $\mathbb U(n)$ (resp. $\mathbb O(n)$). At
last we treat  the permutation process.

\subsection{The one-parameter process}

\subsubsection{Alternate proof of Theorem~\ref{onedim} (annealed)}

The process
\[\mathcal G\sn := \left\{n^{1/2}  \sum_1^n |U_{i1}|^2 \left(1_{R_i \leq s}
  - s\right)\ , \ s\in [0,1]\right\}\] is an example of a so-called weighted
empirical  process.  We could  then  apply  Theorem  1.1 of  Koul  and
Ossiander \cite[p.\ 544]{KoulOssiander}, and use the representation of
$(|U_{11}|^2, \dots, |U_{i1}|^2, \dots, |U_{in}|^2)$ by means of gamma
variables  to check  conditions therein.  But we  prefer to  give (the
sketch of) a proof that is more self-contained and closer 
 to what will
happen in the two-dimensional case.

A possible method for the finite dimensional convergence of $\mathcal G\sn$ is to use the Lindeberg
strategy of replacement  by Gaussian variables. It says  that if $G_1,
\dots,  G_n$ are
 independent  Brownian bridges,  then $\{\mathcal
G\sn (s) \ , \ s\in [0,1]\}$
and $\{n^{1/2} \sum_1^n |U_{i1}|^2 G_i(s)\  , \ s\in [0,1]\}$ have the
same   limits   if    \ben\label{e3.1}\lim_n   \mathbb   E\left[\sum_1^n
  \left(n^{1/2}  |U_{i1}|^2\right)^3\right] =  0\een which  holds true
since $\mathbb E |U_{i1}|^6 = O(n^{-3})$. Then it remains to see that
\[\left\{n^{1/2}  \sum_1^n  |U_{i1}|^2  G_i(s)\   ,  \  s\in  [0,1]\right\}  \law
\left\{\left(n \sum_1^n |U_{i1}|^4\right)^{1/2} G_1(s)\  , \ s\in [0,1]\right\}\]
and to prove \ben
\label{inp}
\lim_n n  \sum_1^n |U_{i1}|^4 =  1 + \beta'^{-1} \een  in probability.
This latter task  may be performed using moments of  order one and two
of the above expression. We skip the details.

To prove  tightness, we  revisit  criterion (14.9)  of Billingsley
\cite{Bill}. 
For $r <s<t$, we have
\ben  \nonumber  \mathbb   E    \left[\left(\mathcal    G\sn(s)   -    \mathcal
    G\sn(r)\right)^2     \left(\mathcal     G\sn(t)     -     \mathcal
    G\sn(s)\right)^2\right]
=O((s-r)(t-s))\\
\times  \mathbb  E \left(\sum_{i\not=j}|U_{i1}|^2|U_{j1}|^2  +  \sum_i
  |U_{i1}|^4\right)  \een  Since $\sum_1^n  |U_{i1}|^2  =  1$ and  $\lim_n
\mathbb E  \left( n \sum_1^n  |U_{i1}|^4\right) = 1  + \beta'^{-1}$,we
have \ben
\label{sup}
\sup_n  \mathbb  E \left(\sum_{i\not=j}|U_{i1}|^2|U_{j1}|^2  +  \sum_i
  |U_{i1}|^4\right) < \infty \een and the proof is ended.

\subsubsection{Conjecture 1}

Inspecting the  above proof, we see  that if all the  convergences and
bounds for  statistics built from $U$  (i.e. (\ref{e3.1}), (\ref{inp})
and  (\ref{sup}))  were  almost  sure,   we  could  claim  a  quenched
convergence  in  distribution.  Of  course, we  probably  need  higher
moments calculus.

\noindent{\bf  Conjecture  1.}  For   $\pi$ - almost  every  $u$,  the
push-forward of $d\omega$ by the mapping
\[\omega \mapsto \sqrt n \left(\mathcal B\sn(\omega, U) - I\right)\]
converges weakly to the distribution of $\sqrt{1+ {\beta'}^{-1}} B_0$.

\subsection{The Haar process}

\subsubsection{Alternate proof of Theorem~\ref{q-a} (annealed)}

For a complete proof in this  flavor, see \cite{arXiv}. We start from
the  following   decomposition,  analogous  to   (\ref{decnew}):  \ben
\mathcal T\sn  - \mathbb E  \mathcal T\sn = \mathcal  V\sn +\widetilde
S\sn\otimes     \widetilde{S'}\sn     +     n^{1/2}\left(     I\otimes
  \widetilde{S'}\sn  +  \widetilde  S\sn\otimes  I\right)  \een  where
$\mathcal V\sn$ is  defined in (\ref{defcalV}), or  explicitly by \ben
\mathcal   V\sn_{s,t}(\omega,  U)   =  \sum_{ij}   \left(|U_{ij}|^2  -
  n^{-1}\right) [\mathbf  1_{R_i \leq s}-s][\mathbf 1_{C_j  \leq t}-t]
\een (compare with $\widehat W\sn$). The scheme consists in proving
\begin{enumerate}
\item  the  convergence  of  $\mathcal  V\sn$  to  $\sqrt{\beta'^{-1}}
  W^{[\infty)}$ in the sense of finite dimensional distributions;
\item  the  tightness  of  the  sequence  $n^{-1/2}\mathcal  V\sn$  in
  $D([0,1]^2)$.
\end{enumerate}
It seems to  be similar to the one above  in the one-dimensional case,
nevertheless  there  are  two  differences. First,  we  do  not  study
$\mathcal T\sn$, but  $\mathcal V\sn$ which is  its ``involved'' part.
Second,  we  did  not  succeed  to prove  directly  the  tightness  of
$\mathcal V\sn$ but only that  of $n^{-1/2}\mathcal V\sn$. Actually we
know  that a stronger result holds true, 
  since  as  a  consequence of  Theorem~\ref{q-a} 
(quenched),  $\mathcal V\sn  \convlaw \sqrt{\beta'^{-1}}W^{(\infty)}$,
but using Theorem~\ref{q-a} here defeats the point.

For the finite dimensional convergence we use again the Lindeberg strategy
replacing first the  processes $1_{R_i \leq s}-s$  by Brownian bridges
$\beta_i(s)$ and afterwards replacing the processes $1_{C_j \leq t}-t$
by Brownian  bridges $\widetilde\beta_j(t)$. The original  process and
the new one have the same limit as soon as \ben
\label{fidi2}
\lim_n  \mathbb E  \sum_{i,j=1}^n |U_{ij}|^6  = 0\,  \een (which  is true,
since again $\mathbb  E |U_{ij}|^6 = O(n^{-3}))$. To  simplify, let us
explain what happens for the  one-dimensional marginal after the above replacement.

 Let $X = (X_1,
\dots,  X_n)$ and  $Y= (Y_1,  \dots, Y_n)$,  two independent  standard
Gaussian  vectors in  $\mathbb R^n$.  We thus  study the  bilinear non
symmetric form
\[Q_n := \sum_{i,j=1}^n X_i (|U_{ij}|^2 - 1/n)Y_j\] built from the non
symmetric matrix  $V= \left(|U_{ij}|^2 - 1/n\right)_{i,j  \leq n}$. 
  The  characteristic  function of $Q_n$  is computed by  conditioning upon $X$ and $U$;  taking  into
account the  Gaussian distribution of $Y$ we are lead  to study
the quadratic form
\[\widehat  Q_n  :=  \sum_{i,j=1}^n  X_iH_{ij}X_j   \  ,  \  H_{ij}  =
(VV^*)_{ij} \ , \ i,j =1, \dots, n\,\]
and  we have to prove   that  \ben\label{cvprobagain}   \lim_n  \widehat   Q_n  =
\beta'^{-1}\,,  \een in  probability. It can be checked with straightforward  calculations of
moments of order 1 and 2 of $\widehat Q_n$, demanding moments of order
8 of the entries of $U$ to be computed.

To  prove  the tightness  of  $n^{-1/2}\mathcal V\sn$  we can  use  a
criterion  of  Davydov and  Zitikis  \cite{DaZ}  (notice that  several
criteria   for   tightness  known   in   the   literature,  such   as
Bickel-Wichura   \cite{bickel1971convergence}  or   Ivanoff  \cite{Iv},
failed    in    this    model).    A    sufficient    condition    is:
\ben\label{last}\mathbb  E \left(\mathcal  V_{s,t}\sn\right)^6 \leq  C
(\max(s,t))^3\een as soon  as $\max(s,t) \geq n^{-1}$.  With a careful
look at dependencies,
we   reach:  \ben\label{to_prove_V}   \mathbb  E^U   \left(\mathcal
  V_{s,t}\sn\right)^6 =  \sum_{i_k, j_k, k=1,  \ldots 6}(\prod_{k=1}^6
V_{i_kj_k}) \E(\prod_{k=1}^6 B_{i_k}) \E(\prod_{k=1}^6 B'_{j_k}).
\een Since the  $B_i$ and $B'_j$ are independent and  centered, in the
RHS of the above equation, the  non-zero terms in the sum are obtained
when the $i_k$ (resp.\ the $j_k$) are equal at least 2 by 2. Using the
following properties:
\begin{itemize}
\item $ |\E((B_i)^k)| \leq \E (B_i)^2 \leq s$ for $2\leq k\leq 6$
\item $ |\E((B'_j)^k)| \leq \E(B'_j)^2 \leq t$ for $2\leq k\leq 6$\,,
\end{itemize}
it was  checked in \cite{arXiv} that (\ref{last}) holds true.

\subsubsection{Conjecture 2}

Inspecting the  above proof, we see  that if all the  convergences and
bounds   for   statistics   built  from   $U$   (i.e.   (\ref{fidi2}),
(\ref{cvprobagain}) and (\ref{last})) were almost sure, we could claim
a  quenched  convergence in  distribution.  This  would require  sharp
analysis of  homogeneous polynomials  in the entries  of $V$  hence of
$U$. \medskip

\noindent{\bf  Conjecture   2.}  For  $\pi$  - almost   every  $u$,  the
push-forward of $d\omega$ by the map
\[\omega \mapsto n^{-1/2}\left(\mathcal T\sn(\omega, U\sn) - nI\otimes I\right)\]
converges weakly  to the  distribution of $\mathcal  W^{(\infty)}$. In
other words $\pi(\mathcal E) = 1$. \medskip

\subsection{The permutation process}

For  the permutation  process  we  have a  complete  picture, i.e.\  a
convergence  conditionally   upon  $u$,   hence  a  direct   proof  of
Theorem~\ref{twodimpermut} (annealed).

\begin{thm}
  For every $u\in \otimes_{n=1}^\infty \mathcal S_n$, the push-forward of $d\omega$ by the mapping
  \[\omega \mapsto n^{-1/2}\left(\mathcal T\sn(\omega, U\sn) - n I\otimes I\right) \]
  converges weakly to the distribution of $B_{00}$.
\end{thm}

\proof If  $\sigma\sn$ is the  permutation associated with  $U\sn$, we
have
\[\mathcal T\sn_{st} = \sum_{i=1}^n {\mathbf 1}_{R_i \leq s}{\mathbf 1}_{C_{\sigma\sn(i)}\leq t}\,.\]
If we fix $u$, we fix
$\sigma\sn$   for  every   $n$.  It   is  clear   that  the   sequence
$C_{\sigma\sn(i)},1\leq i \leq n$ has the same distribution as $C_i, 1
\leq i \leq n$. We have then (conditionally)
\[n^{-1/2}(\mathcal  T\sn (.  , U\sn)  - n  I\otimes I)  \law \mathcal
X_n\] where
\[\mathcal X_n :=
\left(  n^{-1/2}\left(\sum_{i=1}^n  {\mathbf 1}_{R_i  \leq  s}{\mathbf
      1}_{C_i\leq t} -  nst\right) \ , \ s,t \in  [0,1]\right)\,\] is a
classical two-parameter empirical process.  In \cite{Neu} it is proved
that 
this process converges in distribution to $B_{00}$. If $\mathcal F$ is
any bounded continuous  function of $D([0,1)^2)$ in $\mathbb  R$ we may
write, for every $u$,
\[\mathbb E\left[\mathcal F\left(n^{-1/2}(\mathcal T\sn (. , U\sn) -
    n I\otimes  I)\right) \ \middle|  \ u\right] = \mathbb  E \mathcal
F(\mathcal X_n))\rightarrow_{n\rightarrow  \infty} \mathbb  E \mathcal
F(B_{00})\]
which concludes the proof. $\Box$


\begin{thebibliography}{0}

\bibitem{AndFar}
G.W.~Anderson, B.~Farrell.
\newblock{Asymptotic liberating sequences of random unitary matrices}.
\newblock{\em arXiv\/} 1302.5688v2 [mathPR] 2013.


\bibitem{BaCoS}
T.~Banica, B.~Collins and J.M.~Schlenke.
\newblock{On polynomial integrals over the orthogonal group}.
\newblock{\em J. Combin. Theory Ser. A },  118(3)  778-795, 2011.

\bibitem{FBG}
F.~Benaych-Georges.
\newblock{A universality result for the  global fluctuations of the eigenvectors of Wigner matrices}.
\newblock{\em Random Matrices: Theory Appli.\/} 1,   2012, no. 4, 1250011




\bibitem{unistoch}
G.~Berkolaiko.
\newblock{Spectral gap of doubly stochastic matrices generated from equidistributed unitary matrices}.
\newblock {\em J. Phys. A}, 34(22), 319-326,  2001.




\bibitem{bickel1971convergence}
P.J.~Bickel and M.J.~Wichura.
\newblock {Convergence criteria for multiparameter stochastic processes and
  some applications}.
\newblock {\em Ann. Math. Statist.}, 42(5), 1656--1670, 1971.

\bibitem{Bill}
P~Billingsley.
\newblock{Convergence of Probability measures},
\newblock{Wiley} 1968.

\bibitem{AliB}
A.~Boufferoum.
\newblock{Eigenvectors of sample covariance matrix: universality of global fluctuations}
\newblock{\em arXiv: 1306.4277}, 2013.

\bibitem{B2N}
P.~Bourgade, J.~Najnudel and A.~Nikeghbali.
\newblock{A unitary extension of virtual permutations}.
\newblock{\em Int Math Res Notices\/} first published online July 16, 2012 doi:10.1093/imrn/rns167.


\bibitem{chapuy}
G.~Chapuy.
\newblock{Random permutations and their discrepancy process}.
\newblock{\em 2007 Conference on Analysis of Algorithms, AofA 07}
\newblock{\em Discrete Math. Theor. Comput. Sci. Proc. AH\/} 415--426, 2007.


\bibitem{Collins}
B.~Collins.
\newblock{Product of random projections, Jacobi ensembles and universality problems arising from free probability.}
\newblock{\em Probab. Theory Related Fields}, 133, 315--344, 2005.




\bibitem{chaf}
D.~Chafai.
\newblock{The Dirichlet Markov Ensemble},
\newblock{\em J. Multiv. Anal.},  101(3): 555-567, 2010.

\bibitem{DaZ}
Y.~Davydov and R.~Zitikis.
\newblock{On weak convergence of random fields}.
\newblock{\em Ann. Inst. Statist. Math.}, 60, 345--365, 2008.

\bibitem{DPY}
P.~Deheuvels, G.~Peccati and M.~Yor.
\newblock{On quadratic functionals of the Brownian sheet
and related processes}.
\newblock{\em Ann. Probab.}, 116, 493-538, 2006.

\bibitem{CDMAR}
C.~Donati-Martin and A.~Rouault.
\newblock{Truncation of  Haar unitary matrices, traces and bivariate Brownian bridge}.
\newblock{\em Random Matrices Theory Appli.}, vol. 1, no.1, 115007, 24p, 2012.

\bibitem{arXiv}
C.~Donati-Martin and A.~Rouault.
\newblock{Random truncations of  Haar unitary matrices and  bridges}.
\newblock{\em arXiv.} 1302.6539 [mathPR] 2013.

\bibitem{DuPa} I.~Dumitriu and E.~Paquette.
\newblock{Global fluctuations for linear statistics of $\beta$-Jacobi ensembles}.
\newblock{\em Random Matrices Theory and Appli.}, vol. 1, no. 4, 1250013, 60p, 2012.


\bibitem{Fararxiv}
B.~Farrell.
\newblock{Limiting Empirical Singular Value Distribution of Restrictions of Unitary Matrices}
\newblock{\em arXiv\:} 1003.1021v1 2010 [mathFA] 2010.

\bibitem{Farfourier}
B.~Farrell.
\newblock{Limiting empirical singular value distribution of restrictions of discrete Fourier transform matrices}.
\newblock{\em J. Fourier Anal. Appl.}, 17, 733--753, 2011.


\bibitem{Iv}
B.G. Ivanoff.
\newblock {Stopping times and tightness for multiparameter martingales.}
\newblock {\em Statist. Probab. Lett.}, 28, 111--114, 1996.

\bibitem{JiangMax}
T.~Jiang.
\newblock{Maxima of entries of Haar distributed matrices.}
\newblock{\em Probab. Th. Rel. Fields.}, 131, 121-144, 2005.

\bibitem{KoulOssiander}
H.L.~Koul and M.~Ossiander.
\newblock{Weak convergence of randomly weighted dependent residual
              empiricals with applications to autoregression}.
\newblock{\em Ann. Statist.}, vol 22 (1), 540--562, 1994.


\bibitem{Neu}
G.~Neuhaus.
\newblock{On weak convergence of stochastic processes with
              multidimensional time parameter}.
\newblock{\em Ann. Math. Statist.}, 42, 1285--1295, 1971.


\bibitem{Silver1}
J.W.~Silverstein.
\newblock{Describing the behavior of eigenvectors of random matrices using sequences of measures on orthogonal groups}.
\newblock{\em SIAM J.  Math. Anal.} 12(2), 274--281, 1981.


\bibitem{Biao}
B.~Wu.
\newblock{On the weak convergence of subordinated systems}.
\newblock{\em Stat. Probab. Lett.}, 78 3203-3211,  2008.

\end{thebibliography}
\end{document}